\documentclass[11pt, a4paper, twosided]{article}
\usepackage{setspace}
\usepackage{amssymb,amsmath,amsthm,latexsym,bbm,mathrsfs,nicefrac,doc,pdfpages}
\usepackage{a4wide,microtype} 
\usepackage{verbatim} 
\usepackage{hyperref} 
\usepackage{natbib} 
\bibpunct{(}{)}{;}{a}{,}{,} 
\usepackage{graphicx} 
\graphicspath{{./Simulations/}} 

\usepackage[tableposition=top]{caption} 
\newtheorem{theo}{Theorem}[section]

\newtheorem{prop}[theo]{Proposition}
\newtheorem{defi}[theo]{Definition}
\newtheorem{lemm}[theo]{Lemma}

\title{Multi-scale detection of rate changes in spike trains with weak dependencies}
\author{
{\sc Michael Messer$^1$} \and {\sc Kau$\hat{\text{e}}$ M.~Costa$^2$}\and {\sc Jochen Roeper$^2$}\and {\sc Gaby Schneider$^1$}\\[1ex]
{$^1$Institute of Mathematics and $^2$Institute of Neurophysiology}\\{Johann Wolfgang Goethe University}\\
}
\date{}

\begin{document}
\maketitle 
\begin{center}
Accepted for publication in the Journal of Computational Neuroscience,\\
DOI: 10.1007/s10827-016-0635-3. The final publication is available at http://link.springer.com. 
\end{center}
\begin{abstract}
The statistical analysis of neuronal spike trains by models of point processes often relies on the assumption of constant process parameters. However, it is a well-known problem that the parameters of empirical spike trains can be highly variable, such as for example the firing rate. In order to test the null hypothesis of a constant rate and to estimate the change points, a Multiple Filter Test (MFT) and a corresponding algorithm (MFA) have been proposed that can be applied under the assumption of independent inter spike intervals (ISIs).

As empirical spike trains often show weak dependencies in the correlation structure of ISIs, we extend the MFT here to point processes associated with short range dependencies. By specifically estimating serial dependencies in the test statistic, we show that the new MFT can be applied to a variety of empirical firing patterns, including positive and negative serial correlations as well as tonic and bursty firing. The new MFT is applied to a data set of empirical spike trains with serial correlations, and simulations show improved performance against methods that assume independence. In case of positive correlations, our new MFT is necessary to reduce the number of false positives, which can be highly enhanced when falsely assuming independence. For the frequent case of negative correlations, the new MFT shows an improved detection probability of change points and thus, also a higher potential of signal extraction from noisy spike trains. 
\end{abstract}

\section{Introduction}
One fundamental property of neurons is that they can code information by varying their firing rate. In many neurons, information is encoded in transient spike rate changes against a variable baseline. The statistical structure of baseline firing in spike trains is an important determinant of neuronal information encoding \citep{Luczak2013,Hartmann2015}. Likewise, state-related changes in baseline firing rate alone affect signal to noise ratio and the quality of information encoding \citep{Lee2012}. In this context, statistical estimation of rate change points is an important tool for extracting relevant features from neuronal signals, especially during so called spontaneous firing neuronal activity during sleep, periods of quiet wakefulness, under anesthesia, or whenever there is no direct behavioral or sensorial trigger for the recorded neural signal, but which does carry important information about the structure of neuronal networks and the biophysics of individual neurons \citep{Schiemann2012,Luczak2013,Hartmann2015}.

General point process models have been proposed for the description of varying firing rates and their dependence on past spiking activity, external stimuli or behavioral events \citep{BrownKM04,Koyama2010,Pillow2008,Paninski2004,Trucculo2004}. 
The present article is motivated by the observation that in addition to representing potentially important neuronal signals, changes in the firing rate can often have a crucial impact on a large number of standard statistical spike train analyses that require the assumption of a constant firing rate \citep[e.g.,][]{Brody1999,Gruen02b,Schneider08}. Therefore, the main aim is to present a statistical test of the null hypothesis of constant rate and a method that can estimate change points in the firing rate in order to divide a spike train into sections of approximately constant rate.

\begin{figure*}[htbp]
\begin{center}
\includegraphics[width=0.75\textwidth]{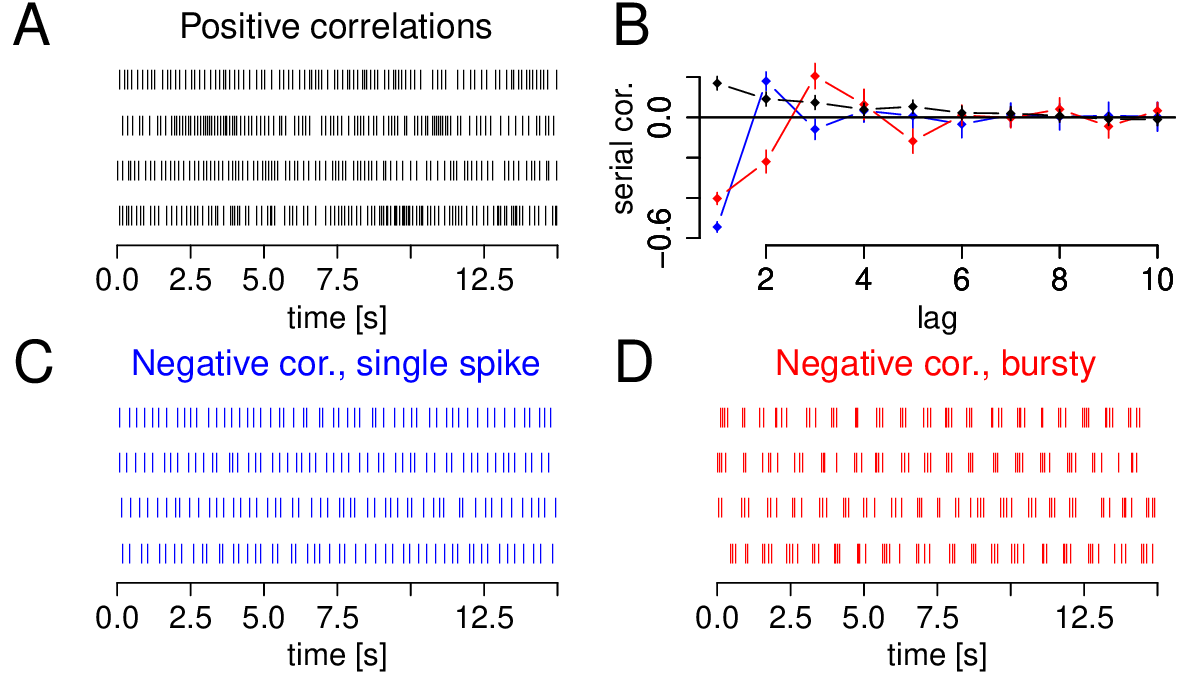}
\caption{Spike trains with serial correlations of ISIs from the sample data set. (A) Positive correlations, (C, D) negative serial correlation of order one. (B) Serial correlations of ISIs, derived from mean in disjoint windows of $50$ ISIs each, with $95$\% confidence limits.}
\label{fig:examples}
\end{center}
\end{figure*}

A statistical method that aims at detecting rate change points in neuronal spike trains should take into account several phenomena and challenges observed in empirical data (see also Figure \ref{fig:examples}). First, distributions of inter spike intervals (ISIs) can be highly diverse, and rate changes can occur on different time scales. Second, other process parameters such as the variance are not known in practice. And finally, as one of the main issues in the present paper, neuronal spike trains have often been reported to show serial dependencies of low orders \citep{Lowen1992,Ratnam2000,Chacron2001,Nawrot2007,Farkhooi2009}, implying that independence of ISIs can not necessarily be assumed in practice. 

Serial correlations themselves have been proposed to be a crucial aspect of information transmission in neuronal spike trains, for example, by reducing variability of spike count through negative correlations, thus increasing signal detection efficiency by a post-synaptic neuron \citep{Chacron2001,Ratnam2000,Chacron2004,Nawrot2007}. Several models of neuronal information coding have been proposed that incorporate mechanisms for positive and/or negative serial ISI correlations \citep{Avila2011,Schwalger2013,Shiau2015}. 
The concept that serial correlations shape the way in which informative spike changes are detected by neuronal systems inspired us to develop a statistical analysis method that can detect rate changes in spike trains while assuming and incorporating serial ISI correlations. Specifically, our novel method can detect rate changes in spike trains with short range dependencies in which the covariance structure of life times is unknown and rate changes may occur at different time scales.

As a proof of principle, we applied our analysis to spike train recordings obtained from spontaneous activity of DA neurons in anaesthetized mice. These trains include more or less regular single spike or bursty patterns, \cite[e.g.,][see Fig.~\ref{fig:examples} C, D for examples]{Bingmer2011,Schiemann2012}, and the activity of this class of neurons has been previously described with spike train models with serial dependencies, such as stochastic cluster processes \citep{Bingmer2011} or Hidden Markov Models \citep{Camproux1996}. The dataset presented here shows serial correlations between successive ISIs, which can be strong for small lags, but decay fast towards zero, in accordance with the literature on serial ISI correlations (Fig.~\ref{fig:examples}).

For the detection of rate changes in  point processes, methods were developed e.g., by \citet{kendall1980,Csorgo1987,Steinebach1993,Gut2002,Gut2009,Messer2014}. 
In the context of stochastic time series, 
change point detection techniques that allow for dependencies have been developed by, e.g., \citet{Tang1993,Lavielle1999,Ray2002,Berkes2006,Dehling2013}. Further, several interesting methods that focus on the aspect of multi-scale detection have been proposed recently by \citet{Munk2014,Fryzlewicz2014,Messer2014}.

Our novel method extends the multiple filter test (MFT) proposed in \cite{Messer2014} to respecting weak dependencies in the ISIs. Regarding the point process model the ISIs are often referred to as the \textit{life times} of the point events. The MFT was designed specifically for spike trains with a wide range of ISI distributions and multi scale rate changes. The  idea in the corresponding filtered derivative approach is to study for every time point the difference between the number of spikes in the left and right window, scaled by an estimate of its standard deviation. Under specific assumptions, one obtains a limit process that is independent of all parameters of the underlying spike train. This limit process can  be used to define rejection thresholds of the null hypothesis of constant rate  and to estimate the rate change points. By simultaneous application of multiple moving windows, change points at multiple time scales can be detected. A corresponding algorithm \citep[MFA, see][]{Messer2014} can then be used to estimate the change points.

The idea behind extending the MFT to weak dependencies is based on the fact that under independence, the variance of the life times $\{\xi_i\}_{i\ge 1}$, $\mathbb Var(\xi_1)$, is used as a scaling factor of the test statistic. If independence does not hold, this term needs to be replaced by 
$$\rho^2:=\mathbb Var(\xi_1)+2\sum_{\ell=1}^\infty\text{Cov}(\xi_1,\xi_{1+\ell}).$$
We then require a consistent estimate of $\rho^2$ in practical application. Here we focus on the practically important case of $m$-dependence, i.e., when Cov$(\xi_1,\xi_{1+\ell})=0$ for all $\ell>m$, with some $m\in \mathbb N$, which yields consistency of the standard estimators for the summands of $\rho^2$.

The paper is organized as follows. We  first review the ideas of the MFT assuming independence and the corresponding MFA for change point detection in section \ref{sect:classicMFT}. In section \ref{sect:depMFT} we  derive a modification that can be applied to spike trains with weak dependencies. Section \ref{sect:PSEex} gives examples of such theoretical processes to illustrate their  correspondence to neuronal spike trains, particularly including also tonic and oscillatory bursty processes. Section \ref{sect:application} uses simulations to discuss estimation principles of $\rho^2$ and $m$ and practical performance of the proposed method including also a recommendation for the choice of the window size. Particularly, we show that disrespecting serial correlations or globally estimating $\rho^2$ or $m$ can yield erroneous results, and illustrate improved performance of the modified MFT and MFA with regard to the number and location of  change points. In section \ref{sect:dataana}, we apply the derived statistical method and algorithms to a data set of  spike train recordings obtained from spontaneous activity of DA neurons in anesthetized mice. 

\section{Extension of the Multiple Filter Test to Weak Dependencies}\label{sect:overview}
We consider a finite spike train of length $T>0$ on the time interval $[0,T]$ as a sequence of spikes
$0 <  S_1 <  S_2 < \cdots<S_{N_T}$, where $N_t$ denotes the number of spikes up to time $t$. The ISIs are denoted by $\{\xi_i\}_{i\ge 1}$, with $\xi_1 = S_1$ and
$$\xi_i = S_i - S_{i-1}\quad \text{for}\quad  i=2,3,\ldots,N_T.$$ 
The ISIs are considered realizations of random variables, and the aim is to construct a statistical test for the null hypothesis that the (positive) mean of all ISIs, i.e., the firing rate, is constant, 
\begin{align}
H_0:\; 	& \mathbb{E}[\xi_i]= \mathbb{E}[\xi_1]=: \mu>0\quad\textrm{for all}\; i=1,\ldots,N_T.
\end{align}	
For the alternative of $k$ change points $c_1,\ldots, c_k\in[0,T]$, we assume $k+1$ (independent) processes with constant rates $\mu_1^{-1},\ldots,\mu_{k+1}^{-1}$, while $\mu_j \not=\mu_{j+1}$ for all $j$. At time zero start in first process with rate $\mu_1^{-1}$, at the first change points $c_1$ jump into the second process of rate $\mu_2^{-1}$ etc. Then the resulting process is a piecewise combination of sections with different rates. If the null hypothesis is rejected, we are interested in estimating the change points $c_1,\ldots,c_k$ in order to segment the spike train into sections of constant rate. 

\subsection{The MFT for rate changes in renewal processes}\label{sect:classicMFT}
Here we describe the main idea of the MFT \citep[for more details see][]{Messer2014}.
The MFT is based on a filtered derivative approach that compares the numbers of events, $N_{\text{le}}:=N_{t}-N_{t-h}$ and $N_{\text{ri}}:=N_{t+h}-N_{t}$ in the left and right window of  size $h\in(0,T/2]$ for every time $t\in[h,T-h]$. By standardizing with a consistent estimator of the standard deviation of this difference, $\hat  s_{h,t}$, one obtains a filtered derivative process 
\begin{align}\label{eq:def_ght}
G_{h,t}:= \frac{N_{\text{ri}}-N_{\text{le}}}{\hat s_{h,t}}.
\end{align}	
Large differences between the numbers of events in the left and right window, i.e., large deviations of $G$ from zero indicate deviations from the null hypothesis of constant rate. In order to test statistical significance of these deviations, the maximal deviation $\max_{t} |G_{h,t}|$ from zero could serve as a test statistic for one window, and the rejection threshold at level $\alpha$ can be derived from the limit process of $G$ as follows. Using an extension $G^{(n)}:=G_{nh,nt}$ in an asymptotic setting in which the window size $nh$ and the time $nT$ (or alternatively, the firing rate) grow linearly in $n$, $G^{(n)}$ can be shown to converge weakly to a functional $L$ of a standard Brownian motion $W$, 
\begin{align}\label{limit_result}
L_{h,t}:= \left((W_{t+h}-W_t)-(W_t-W_{t-h})\right)/\sqrt{2h},
\end{align}
under the null hypothesis of a constant rate.
Note that $L_{h,t}\sim N(0,1)$ for all $h$ and $t$, i.e., $\hat s_{h,t}$ standardizes the difference of the number of events in both windows.
As $L$ does not depend on parameters of the underlying process, the distribution of $\max_t |L_{h,t}|$ can be easily simulated to obtain a rejection threshold $Q$ for a statistical test at level $\alpha$. 

In order to allow detection of change points at multiple time scales, the MFT combines multiple windows from a finite set $H$ and the corresponding processes $G_{h,t}$. As the distribution of $\max_t |G_{h,t}|$ depends on the window size $h$, the process $G$ is rescaled to give about the same weight to every window size, resulting in a rescaled process
$$R_{h,t}:=\frac{|G_{h,t}|-\hat \mu_{M_h^*}}{\hat \sigma_{M_h^*}},$$
where $\hat \mu_{M_h^*}$ and $\hat \sigma_{M_h^*}$ denote the estimated mean and standard deviation of $M_h^*:=\max_t |L_{h,t}|$ obtained in simulations by simulating $W$ and deriving $L$ from $W$ as in (\ref{limit_result}). The maximum of all $R$-processes,
$$M:=\max_{h,t} R_{h,t},$$
is used as a test statistic. The rejection threshold $Q$ can then be derived from the corresponding distribution of $\max_h{{(M_h^*-\hat \mu_{M_h^*})}/{\hat \sigma_{M_h^*}}}$, which can be obtained in simulations.

This approach has three practical advantages: First, it does not require previous knowledge of process parameters because $G$ is scaled such that the limit process does not depend on the parameters of the underlying process. Second, it can be applied to a wide range of processes, i.e., Poisson or Gamma processes or processes with complex or unknown ISI distributions as long as ISIs are independent and identically distributed \citep{Steinebach1995}. It even holds for processes with independent  but not necessarily identically distributed ISIs in the sense that the variance of ISIs may show a certain degree of variation between regular and irregular phases \citep{Messer2014} as can sometimes be observed in empirical spike trains. Third, this approach allows the simultaneous use of multiple windows in a finite  set $H$ and thus, analysis of change points at multiple time scales. Due to the asymptotic nature of the method, the smallest window should contain at least about $100-200$ spikes in order to approximately keep the significance level.

\begin{figure}[htbp]
\begin{center}
\includegraphics[width=0.45\textwidth]{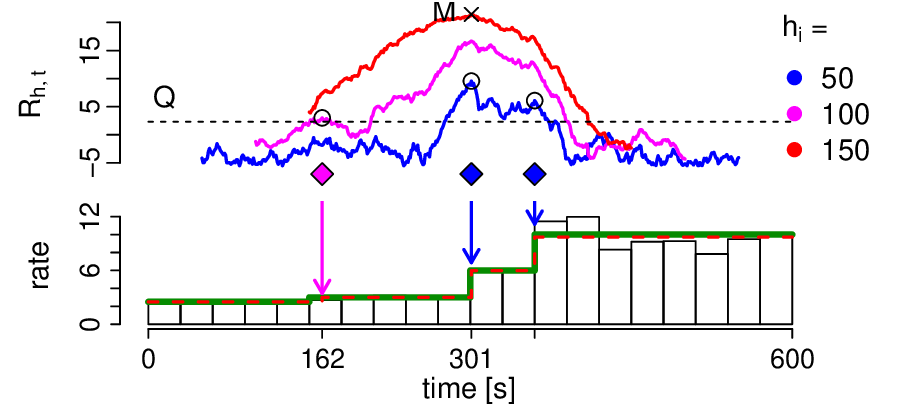}
\caption{Application of the MFT and MFA to a piecewise renewal process with Gamma-distributed intervals $\xi_i$ with variance $\text{Var}(\xi_i)=0.2^2$ and rates $2.5, 3, 6, 10$ Hz with change points at $150, 300, 360$ seconds. Upper panel indicates rescaled processes $R_{h,t}$ for window sizes $h \in  \{50, 100, 150\}$ and rejection threshold $Q$ (dashed) at $5\%$ level, and lower panel indicates rate histogram, true rate (thick, solid) and estimated rate (thin, dashed) with estimated change points (diamonds and arrows). The color indicates the window with which the change point was detected.}
\label{fig:MFT1}
\end{center}
\end{figure}

The MFT is applied to a simulated spike train with three rate changes in Figure \ref{fig:MFT1}. The upper panel indicates the rescaled processes for a window set $H=\{50,100,200\}$. The maximum $M$  exceeds the rejection threshold $Q$, and the null hypothesis of constant rate is rejected. Then, the MFA successively estimates the change points. For every window $h$, change point candidates $\hat c_j$ are  identified by successively locating the maxima of $(R_{h,t})_t$ and then deleting their $h$-neighborhood $[\hat c_j-h,\hat c_j+h]$. Change point candidates are then successively combined (see also the articles by \cite{Fryzlewicz2014,Munk2014} for similar approaches), preferring candidates of smaller windows and adding only those whose $h$-neighborhood does not overlap accepted change points. This is motivated by the idea that large windows tend to be affected by multiple change points, which may reduce their estimation precision. In Figure \ref{fig:MFT1}, change points with fast, strong changes are estimated with small windows, while change points with slow and weak  changes are estimated with larger windows. 

\subsection{The MFT for weak dependencies}\label{sect:depMFT}
The main purpose of this paper is to study the MFT in case of weak dependence of ISIs. We will show here that this requires two assumptions: First, a generalized class of point processes that also include weak dependencies (Definition \ref{def:rpsel}), and second, consistent estimation of process parameters (Proposition \ref{prop:main}).

Under independence, a consistent estimator of the standard deviation of $(N_{\text{ri}}-N_{\text{le}})$ in equation (\ref{eq:def_ght}) is given by
\begin{equation}
\left( 2 n h \hat\sigma^2 / \hat\mu^3 \right) ^{1/2}
\label{eq:sindep}
\end{equation}
if $\hat \mu > 0$, and zero otherwise, where $\hat \mu$ and $\hat \sigma^2$ denote the empirical mean and  variance of the ISI lengths in the analysis window. In case of non-zero covariances between successive ISIs, the estimate of the variance of ISIs $\sigma^2:=\mathbb{V}\!ar(\xi_1)$ needs to be replaced by an estimator of
\begin{equation}\label{eq:rho}
\rho^2 =\sigma^2 + 2 \sum_{\ell=1}^{\infty} \rho_\ell,
\end{equation}
where $\rho_\ell:=\text{Cov}(\xi_1,\xi_{1+\ell})$, yielding 
\begin{equation}\label{hats}
\hat s := \hat s_{nh,nt}:=\left( 2 n h \hat\rho^2 / \hat\mu^3 \right)^{1/2},
\end{equation}
where details on $\hat \rho^2$ and $\hat \mu$ can be found in section \ref{sect:consistentest}.
Here we show under which assumptions on the point processes and modifications of the MFT one obtains the same convergence and thus, applicability to spike trains with weakly dependent ISIs. To that end, we require a class of point processes $\mathscr P$ for which the ISIs fulfill a functional central limit theorem (FCLT) and for which consistency of $\hat s$ can be concluded. 

\begin{defi}{}\label{def:rpsel}\quad\\
The class of point processes $\mathscr P$ is given by all point processes on the positive line whose life times $\{\xi_i\}_{i\ge 1}$ are stationary, ergodic, almost surely positive and square-integrable and further they fulfill $\rho^2 >0$ (see (\ref{eq:rho})) as well as 
\begin{align}\label{sumcond}
\sum_{i=2}^{\infty} \|\mathbb{E}[\xi_1 -\mathbb{E}[\xi_1]|\{\xi_k|k\ge i\}] \| <\infty.
\end{align} 
\end{defi}

Here, $\|\cdot\|$ denotes the $L^2$-norm (as the conditional expectation in (\ref{sumcond}) is a random variable). Stationarity means that the distribution of any subset of life times is invariant under a time shift of their indices. The assumptions on  $\{\xi_i\}_{i\ge 1}$  particularly imply a FCLT as well as ergodic theorems. The FCLT will be used to derive the convergence of the filtered derivative process (Proposition \ref{prop:main}) and the ergodic theorems will be used for consistent parameter estimation (Lemma \ref{mdependent1} and \ref{mdependent2}). See \cite{Billingsley1999} for details on the notions of stationarity and ergodicity. Further note that the summation condition (\ref{sumcond}) implies absolute convergence of the series (\ref{eq:rho}) \citep[see][Thm.~19.1]{Billingsley1999}. This condition particularly holds true for the special case of $m$-dependent sequences.

Throughout this article, $\stackrel{d}{\longrightarrow}$ denotes convergence in distribution, and $(D[0,\infty),d_{SK})$ denotes the space of c\`adl\`ag-functions on $[0,\infty)$ endowed with Skorokhod-topology, and analogous for $(D[h,T-h],d_{SK})$. 

The following proposition ensures that the MFT can be applied to point processes $\Phi\in\mathscr P$ when their parameters are consistently estimated. 
\begin{prop}\label{prop:main}
Let $\Phi\in \mathscr P$ with ISIs $\{\xi_i\}_{i\ge1}$ and let $\hat s^2$ be an estimator for $s^2=2nh\rho^2/\mu^3$ that satisfies in $(D[h,T-h],d_{SK})$ as $n\to\infty$ that $(\hat s/s)_t \to (1)_t$ in probability.\\
Then it holds for the filtered derivative process $G^{(n)}=(G_{nh,nt})_t$ as given in (\ref{eq:def_ght}) in $(D[h,T-h],d_{SK})$ as $n\to\infty$
\begin{align*}
G^{(n)} \stackrel{d}{\longrightarrow} L.
\end{align*}
\end{prop}	

\noindent\textbf{Proof:} 
From the conditions on $\mathscr P$ it follows that in $(D[0,\infty),d_{SK})$ as $n \to \infty$
\begin{align}\label{waiting_time_condition}
\left( \frac{1}{\rho\sqrt{n}} \sum_{i=1}^{[n t]}(\xi_i - \mu)\right)_{t} \stackrel{d}{\longrightarrow}W,
\end{align}
where $W$ denotes a standard Brownian motion \citep[see][Thm.~19.1.]{Billingsley1999}. For $t\ge 0$ let $$Z_t^{(n)} := (N_{n t} - nt/\mu)/(n\rho^2/\mu^3)^{1/2}$$ denote the rescaled counting process. According to \cite{Vervaat1972} it follows from (\ref{waiting_time_condition}) that in $(D[0,\infty), d_{SK})$ as $n\to\infty$ it holds $Z^{(n)} \stackrel{d}{\longrightarrow}W$. We then define a continuous map $\varphi_h: (D[0,\infty),d_{SK}) \to (D[h,T-h],d_{SK})$ via 
$f(t)\stackrel{\varphi_h}{\longmapsto}  ((f(t+h)-f(t)) - (f(t)-f(t-h)))/(2h)^{1/2}$. By continuous mapping theorem it follows in $(D[h,T-h],d_{SK})$ for $n\to\infty$ that 
$$((N_{\text{ri}}^{(n)}-N_{\text{le}}^{(n)})/(2nh\rho^2/\mu^3)^{1/2})_{t}\stackrel{d}{\longrightarrow} L,$$
where $N_{\text{ri}}^{(n)}:=N_{n(t+h)}-N_{nt}$ and $N_{\text{le}}^{(n)}:=N_{nt}-N_{n(t-h)}$. 
Due to the consistency assumption of the estimator $\hat s$, we can exchange  $(2nh\rho^2/\mu^3)^{1/2}$ with $\hat s$ by Slutsky's theorem.\hfill$\Box$\\

\subsection{Examples for practical application}
Proposition \ref{prop:main} states that the MFT is applicable to processes in the class $\mathscr P$ if one uses the modified filtered derivative process
\begin{equation}\label{eq:Greal}
G_{h,t}:=\frac{N_{\text{ri}}-N_{\text{le}}}{\hat s_{h,t}},
\end{equation}
with $\hat s^2$ a consistent estimator of $s^2 = 2h{\rho^2/\mu^3}$ and $\rho^2 = \sigma^2 + 2 \sum_{\ell=1}^{\infty} \rho_\ell$ with the convention $G_{h,t}:=0$ if $\hat s_{h,t}=0$. In order to illustrate practical applicability specifically to spike trains with weakly dependent life times, we give examples of processes in $\mathscr P$ that resemble empirical spike trains  (section \ref{sect:PSEex}) and examples of consistent estimators of $s$ (section \ref{sect:consistentest}).

\subsubsection{Processes in $\mathscr P$}\label{sect:PSEex}
The assumptions of processes in $\mathscr P$ are fulfilled for example in renewal processes with independent and identically distributed ISIs. Here, we focus on dependencies in the ISI structure, i.e., processes with stationary and ergodic ISIs as stated in Definition \ref{def:rpsel}. In a simple but practically important case, the ISIs are $m$-dependent for an $m\in \mathbb N$, i.e., $\rho_{\ell}=0$ for all $\ell>m$.

\begin{figure*}[h!]
\begin{center}
\includegraphics[width=0.75\textwidth]{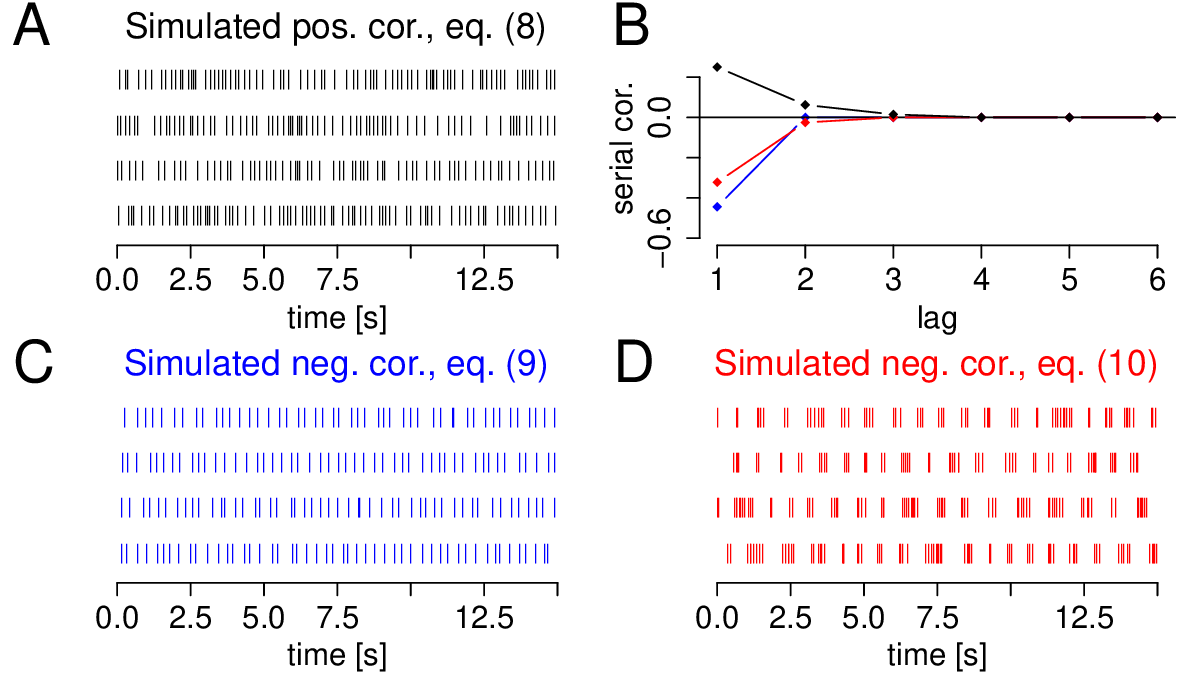}
\caption{Examples of processes from $\mathscr P$ (panels A,C,D) that resemble spike train patterns from Fig.~\ref{fig:examples}, and corresponding covariance structure (B). A. A $3$-dependent process according to equation (\ref{eq:MAmodel}), with $X_i$ Gamma distributed and $\mu=0.2, \sigma=0.1, a_k=c^k, c=0.25$. C. A $1$-dependent process similar to Fig.~\ref{fig:examples} C, simulated according to equation (\ref{eq:ss}), with $\mu=0.3,\sigma_1=0.06,\sigma_2=0.12$. D. A $2$-dependent process similar to Fig.~\ref{fig:examples} D, simulated according to model (\ref{eq:cl}), with $p_I=.5, p_J=.4$, $X_i\sim U[0.45,0.73], Y_i\sim U[.01,.12]$.}
\label{fig:PSEex}
\end{center}
\end{figure*}

Here we give three examples of $m$-dependent processes in Figure \ref{fig:PSEex} that resemble the neuronal spike trains shown in Figure \ref{fig:examples}. Panel A shows a process with $m=3$ and positive serial correlations given by life times
\begin{equation}
\xi_i:= a_0 X_i + a_{1} X_{i-1}+\ldots + a_m X_{i-m},
\label{eq:MAmodel}
\end{equation}
with $X_1,X_2,\ldots$ independent with expectation $\mu_X$ and variance $\sigma_X^2>0$. This implies $$\sigma^2=\text{Var}(\xi_i)= \sigma_X^2\sum_{j=0}^m a_j^2$$ and $\rho_\ell=\sigma_X^2\sum_{j=0}^{m-\ell} a_ja_{j+\ell}$ for $\ell\le m$, and $\rho_\ell=0$ for $\ell>m$ i.e., 
\begin{equation*}
\rho^2 = \sigma_X^2 \left(\sum_{j=0}^m a_j^2 + 2 \sum_{\ell=1}^m \sum_{j=0}^{m-\ell}a_ja_{j+\ell}\right).
\end{equation*}
Appropriate conditions on the $a_i$ and the distribution of $X_i$ ensure almost surely positive ISIs.

Panel C shows an example of a single spike process similar to Figure \ref{fig:examples} C. It is described by ISIs 
\begin{equation} \label{eq:ss}
\xi_i = U_i + Z_i - Z_{i-1},
\end{equation}
where $U_i, Z_i$ are independent and uniformly distributed with $U_i\sim U[\nu-\sigma_1,\nu+\sigma_1]$ and $Z_i\sim U[-\sigma_2,\sigma_2]$, with $\nu, \sigma_1, \sigma_2>0$ and $\sigma_1+2\sigma_2\le\mu$, which assures $\xi_i>0$. In this process, all ISIs $\xi_i$ are identically distributed with mean $\nu$ and variance $\sigma^2=(1/3) (\sigma_1^2+2\sigma_2^2)$. The process is $1$-dependent with negative covariance of lag one given by $\rho_1=-\text{Var}(Z_i)=-(1/3) \sigma_2^2$ (panel B). The spikes of this process can be regarded as jittered uniformly with jitter $Z_i$ around the unobservable beats of a background rhythm with period $\nu$ which is a renewal process with independent and uniformly distributed intervals $U_i$. Related doubly stochastic Cox processes have been used earlier for the description of single spike processes \citep{Bingmer2011}. Similar to Hidden Markov Models \citep{Camproux1996}, they can also be used for the description of oscillatory bursty activity as in Figure \ref{fig:examples} D. 

In order to illustrate applicability of Proposition \ref{prop:main} also to oscillatory bursty spike trains, Figure \ref{fig:PSEex} D shows an example of a $2$-dependent oscillatory bursty process similar to the spike train in Figure \ref{fig:examples} D. Every ISI $\xi_i$ is described by 
\begin{align}\label{eq:cl}
\xi_i & = I_i(1-I_{i-1}) X_i + I_{i-1}J_i Y_i + I_{i-2} J'_i Y'_i \\
&+ (1-\max(I_i(1-I_{i-1}),I_{i-1},I_{i-2} J_i)) Y''_i, \nonumber
\end{align}
where $(I_i)_{i\ge 1}$, $(J_i)_{i \ge 1},(J'_i)_{i \ge 1}$ are independent sequences of independent $\{0,1\}-$valued random variables with success probabilities $p_I$ and $p_J=p_{J'}$, and $(X_i)_{i\ge 1}$, $(Y_i)_{i\ge 1}$,$(Y'_i)_{i\ge 1}$ and $(Y''_i)_{i\ge 1}$ are independent sequences of independent and almost surely positive random variables and $Y_i, Y'_i, Y''_i$ are identically distributed for all $i$. Obviously, all ISIs are identically distributed and the process is $2$-dependent. The idea is that $X_i$ takes large values to generate large ISIs, while $Y_i,Y'_i,Y''_i$ take small values. Then, an ISI  $\xi_i$ takes a large value if $I_i=1$ and $I_{i-1}=0$, such that in this example, a long ISI is typically followed by at least one short ISI, leading to negative serial correlation (panel B). The last summand in (\ref{eq:cl}) only ensures that $\xi_i>0$.

\subsubsection{Consistent estimators}\label{sect:consistentest}
In addition to requiring a process in $\mathscr P$, the second ingredient of Proposition \ref{prop:main} is a consistent estimator $\hat s$. 
Common approaches in the setting of dependencies include methods based on covariance kernel 
estimation \citep[e.g.,][]{deJong2000,Wied2012} or the Bartlett-estimator \citep{Berkes2005,Xiao2012,MuhsalKirch2014}. 
Here we focus on two simple estimators - a \textit{global} and a \textit{local} estimator - that are particularly useful in practical application. Under $m$-dependence, we show consistency under the null hypothesis of constant rate, even for the local estimator. The local estimator is particularly useful in the presence of change points because the global estimator is sensitive to rate changes and therefore tends to be biased in these cases. In contrast, the local estimator does not tend to be biased on most time sections (see section \ref{sect:application}).

In case of $m$-dependence $\rho^2$ equals a finite sum
\begin{align}\label{rho_m}
\rho^2 = \sigma^2 + 2 \sum_{\ell=1}^{m} \rho_\ell. 
\end{align}
The global estimator uses global estimates of the variance and covariances  in (\ref{rho_m})  from the entire spike train using standard estimators
\begin{align}
\hat\rho_\ell &:= \left(\frac{1}{N_{nT}-(\ell+1)}\sum_{i=1}^{N_{nT}-(\ell+1)} \xi_i\xi_{i+\ell}\right) - \hat\mu^2,\label{eq:rhok}\\
\hat\rho^2 & :=\hat\sigma^2+2\sum_{\ell=1}^{m}\hat\rho_\ell,
\end{align}
where $\hat \mu$ denotes the empirical mean of all ISIs. Lemma \ref{mdependent1} in the Appendix shows that this yields a consistent estimator 
\begin{align}\label{s1_mdependent}
\hat s^2 & :=2hn\hat\rho^2/\hat\mu^3
\end{align}
under the null hypothesis. 

As mentioned above, one main disadvantage of  global parameter estimation is that it tends to be biased under the alternative hypothesis (see section \ref{sect:application} and Figure \ref{fig:mestimation} D). Therefore, we suggest to use an analogous local estimator, which for every $t$ uses only the ISIs in the window $(n(t-h),n(t+h)]$. More precisely, for every time  $t$, we estimate $\rho^2$ and $\mu$ analogously, 
but only from the life times that lie within the windows (separate estimation for the left and the right window)
 and let the local estimator be
\begin{align}\label{s2_mdependent}
\hat s^2:=\left(\frac{\hat\rho_{\text{ri}}^2}{\hat\mu_{\text{ri}}^3}+\frac{\hat\rho_{\text{le}}^2}{\hat\mu_{\text{le}}^3}\right)nh.
\end{align}

For the case of independent life times, i.e., $m=0$, consistency of this estimator was shown in \cite{Messer2014}. In Lemma \ref{mdependent2} in the Appendix we show consistency  of this estimator for $m$-dependent processes.

\section{Practical application of the  MFT for weak dependencies}\label{sect:application}
Section \ref{sect:overview} presented theoretically a class of processes, estimators and statistics that allow to apply the MFT and MFA for the estimation of rate change points in spike trains with weakly dependent ISIs. Here we use simulations to illustrate the difference between the proposed method and the classical MFT that assumes independence. In addition, we discuss the important practical issue of estimating serial dependencies and of choosing the set of windows $H$, particularly the smallest window. Simulations are performed using models (\ref{eq:MAmodel}) and (\ref{eq:ss}), which yield flexible and simple formulas for serial correlations. 

For ease of notation we  denote the MFT and MFA that assume $m$-dependence by MFT$^{(m)}$ and MFA$^{(m)}$. The classical MFT assuming independence will therefore be denoted by MFT$^{(0)}$. All procedures use the statistic described in equation (\ref{eq:Greal}). Under $m$-dependence, $\rho^2$ is estimated up to the $m$-th summand in the MFT$^{(m)}$. The corresponding estimator of $s$ is denoted by $\hat s^{(m)}$.

\begin{figure*}[htbp]
	\begin{center}
		\includegraphics[width=0.375\textwidth]{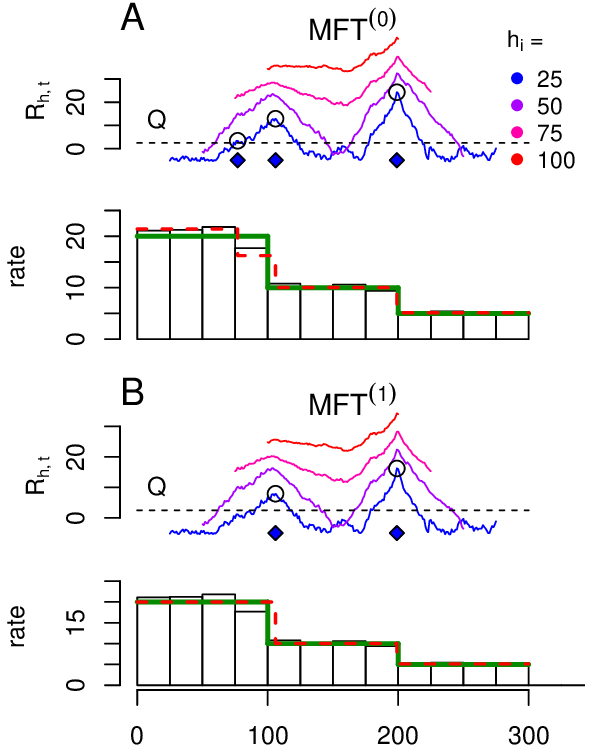}
		\includegraphics[width=0.375\textwidth]{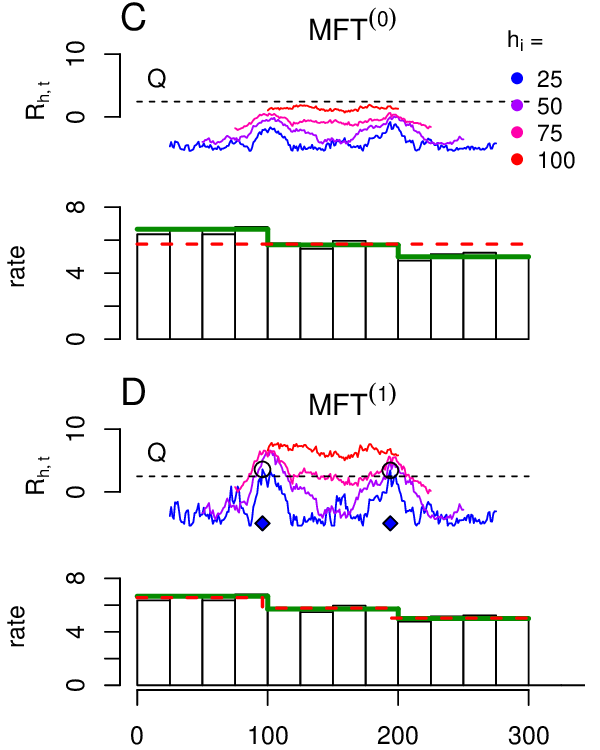}\\[.5ex]
		\includegraphics[width=0.55\textwidth]{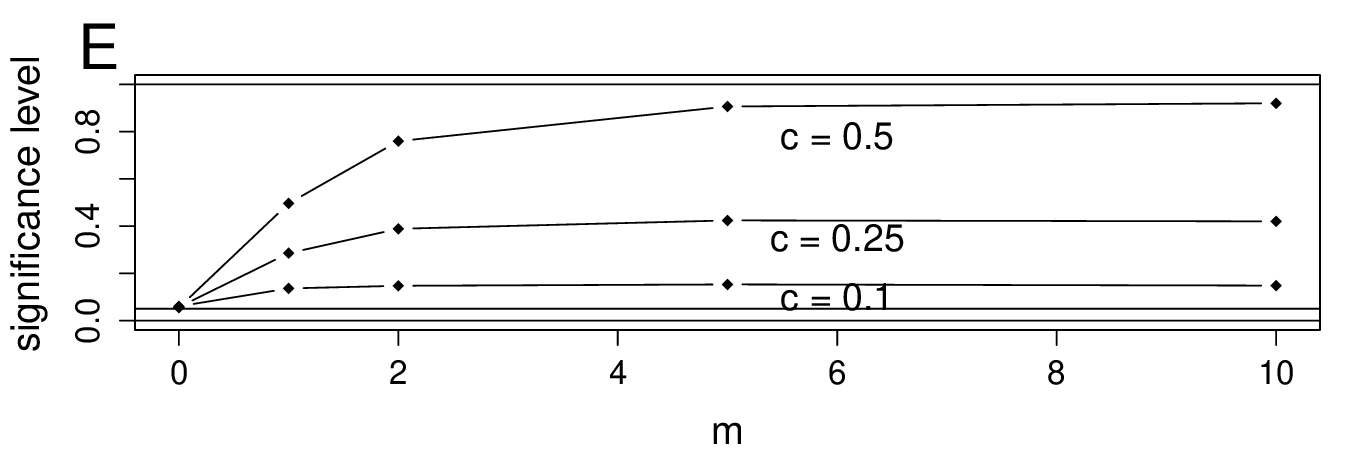}
		\caption{The classical MFT$^{(0)}$ can estimate too many change points when applied to processes with positive correlations (A,E) and may show reduced test power when applied to negative correlations (C). True rate profiles indicated in thick solid, estimated rate profiles in thin, dashed.  Here, the MFT$^{{(1)}}$ (B,D) detects all true change points and no false positives. E. Significance level of classical MFT$^{(0)}$ for positive serial correlations obtained in $10000$ simulations, where the threshold was chosen such that under independence, the MFT$^{(0)}$ would yield an asymptotic significance level of $5\%$ (indicated by horizontal line). All simulations were performed using model (\ref{eq:MAmodel}) with $T=300, H=\{25,50,75,100\}$ and $m=1$ in (A-D) and varying $m$ in E. The coefficients $a_i$ were $a_0=1$ throughout and $a_1=0.5$ in A, B, $a_1=-0.5$ in C,D and $a_k=c^k,c\in\{0.1,0.25,0.5\}$ in E. The $X_i$ were Gamma-distributed in A,B,E and, in order to ensure a.s.~positive ISIs, uniformly distributed in C,D. The parameters $\mu_X,\sigma_X$ were chosen such that the ISIs $\xi_i$ had standard deviation $0.15$ and the given rate profiles (A-D) or $\mu=0.1$ (E). }
		\label{fig:simbsp}
	\end{center}
\end{figure*}

\begin{figure}[htbp]
\begin{center}
\includegraphics[width=0.45\textwidth]{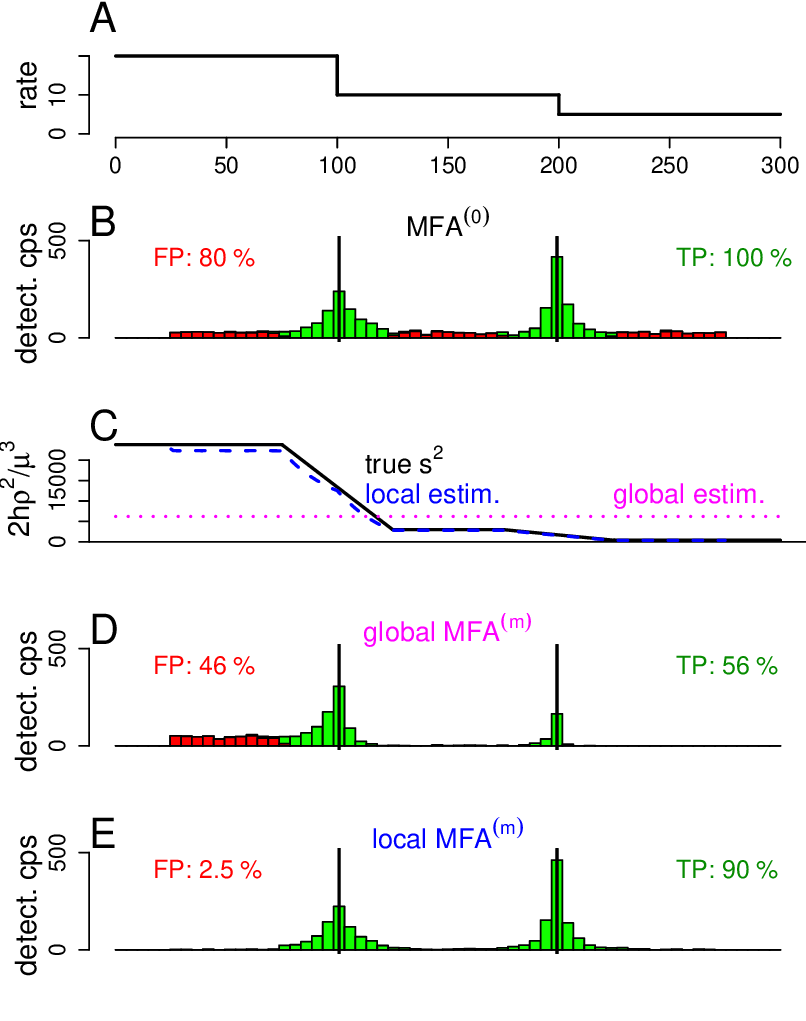}
\caption{Number of detected change points in $1000$ simulations under the alternative with rate profile given in A. B. MFA$^{(0)}$ without accounting for serial correlations. C. True value of local $s^2$ (solid) compared with local estimation (dashed) and global estimation (dotted). D. MFA$^{(m)}$ with global estimate of $s$. E. MFA$^{(m)}$ with estimates of $s$ derived separately in every analysis window. 
Simulations according to model (\ref{eq:MAmodel}), with $X_i$  Gamma distributed, $a_k=0.5^k$, $T=300, H=\{25,50,75,100\},\sigma=0.15, m=3, c=0.5$.}
\label{fig:preferlocal}
\end{center}
\end{figure}

First, we show that falsely applying the classical  MFT$^{(0)}$ yields too many false positives in cases of positive correlations and reduced test power for negative correlations. This is because MFT$^{(0)}$ uses $\hat\rho^2:=\hat\sigma^2$, disrespecting potential serial correlations. Positive correlations yield $\rho^2>\sigma^2$ and  thus increase the number of false positives in the MFT$^{(0)}$ when the scaling $\hat s^{(0)}$ is spuriously low (Figures \ref{fig:simbsp} A, E and \ref{fig:preferlocal} B). Vice versa, negative correlations yield conservative results for the MFT$^{(0)}$ and a reduced test power (Figure \ref{fig:simbsp} C), while in the given example, the MFT $^{(m)}$ can detect the given change points with high precision (Figure \ref{fig:simbsp} D). 

\begin{figure}[htbp]
\begin{center}
\includegraphics[width=0.45\textwidth]{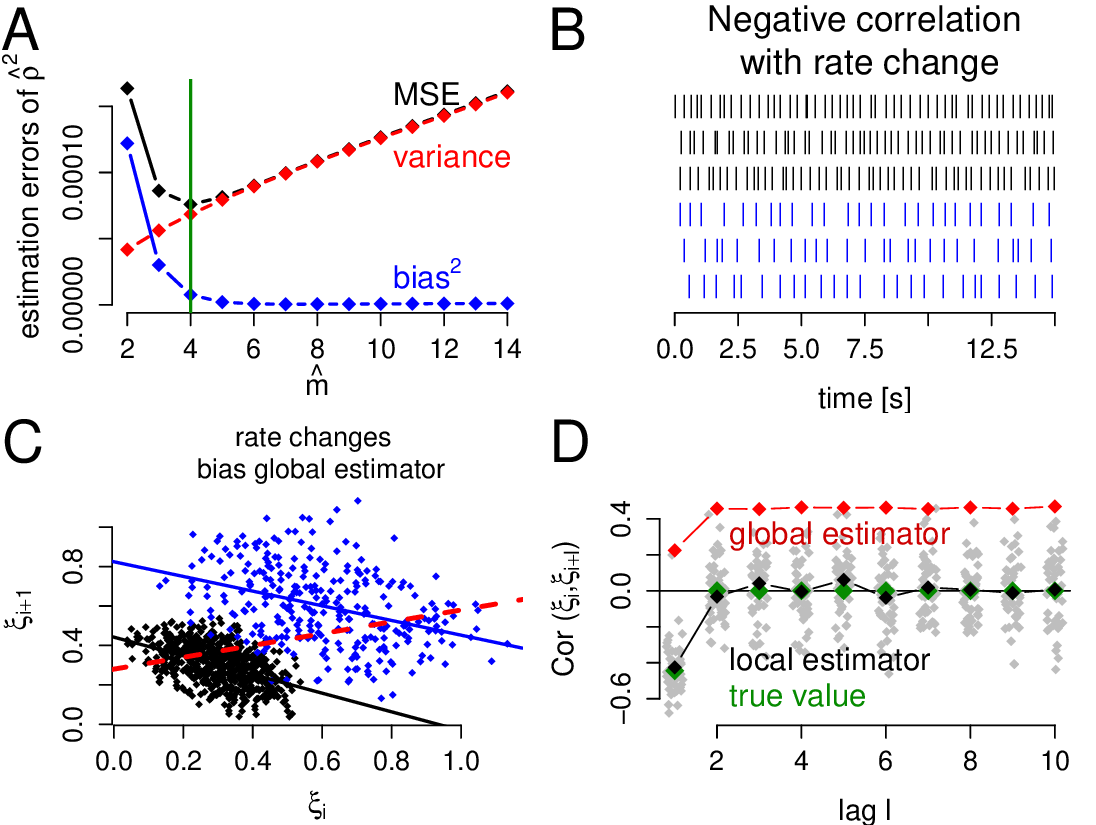}
\caption{A. Estimation errors of $\hat \rho^2$ for different values of $m$. Simulations according to model (\ref{eq:MAmodel}), with $X_i$  Gamma distributed, $a_k=0.5^k$, $T=300,\mu=0.15, \sigma=0.15, m=7$. Serial correlation decreases exponentially with the lag, such that the summands with large lags show only small contributions to $\rho^2$. Neglecting these in the estimation of $\rho$ yields a small bias but a highly reduced variance.
B. A spike train with negative first order serial correlation but a rate change. C. First order correlation is negative locally within sections of constant rate (dots and solid regression lines), but positive globally across sections (dashed regression line). D. Estimates of serial correlations of lags $\ell\in\{1,2,\ldots,10\}$ for the spiketrain shown partly in $B$. Global estimates (red), local estimates (grey) derived in disjoint windows of length $50$ ISIs, their medians (black) and true serial correlations (green).}
\label{fig:mestimation}
\end{center}
\end{figure}

Second, we illustrate the performance of the MFA$^{(m)}$ when $m$ is known using the standard estimators $\hat s^{(m)}$ from section \ref{sect:depMFT} and emphasize that $s$ should be estimated locally. In particular, we propose to use the local estimator from equation (\ref{s2_mdependent}) because a global estimator (eq.~(\ref{s1_mdependent})) would be biased in case of rate changes and thus, reduce test power and/or increase the number of false positives. This effect is illustrated in Figure \ref{fig:preferlocal}. Spike trains with positive serial correlations are simulated with a rate profile with two change points (panel A). As described above, the classical MFA$^{(0)}$ assumes independence and therefore shows many false positives (panel B). Using the MFA$^{(m)}$ with global estimation of $s$ is also unsatisfactory as is shows increased false positive rate on the left and decreased detection rate on the right (panel D). This is because the rate changes cause the true value of $s^2$ to change across time (panel C). The global estimate (dotted) falsely uses a global $\mu$ and therefore a biased global estimate of $\rho^2$  (see also Figure \ref{fig:mestimation} D) and thus underestimates $s^2$ on the left and overestimates on the right. In contrast, the estimates from local windows (blue) correspond closely to the true value of $s^2$, and accordingly, the corresponding MFA$^{(m)}$ using local estimators detects the change points with high precision without showing an increased false positive rate (panel E). For individual examples with positively or negatively correlated life times in the case of $1$-dependence see also Figure \ref{fig:simbsp} B and D). 

Third, we discuss the estimation of $m$, which is typically unknown in practice. If a spike train was arbitrarily long, we could simply use all serial correlations up to an arbitrarily large order as $\rho_\ell=0$ for $\ell>m$, which does not bias the estimation of $\rho$ \citep[this is the idea behind  approaches  for consistent estimation under long-range dependence, see][]{deJong2000,Berkes2005,Wied2012,Xiao2012}. However, in practice, this approach is not applicable because for finite spike trains it highly increases the variance of $\hat \rho$ and thus, the probability of over- or underestimating $\rho$, whereas the former decreases the test power and the latter increases the number of false positives. Therefore, it is  important to include only the largest summands into the estimation of $\rho^2$, while summands with smaller contributions can be neglected. This effect reduces the mean squared error (MSE) of $\hat \rho^2$ by introducing small bias but reducing variance as shown in Figure \ref{fig:mestimation} A where for $m=7$, $\hat m=4$ yields the smallest MSE. 
 
Therefore, we consider here only the practically important case in which serial correlations decrease with the lag, and propose to search the smallest lag $\ell^*$ for which the serial correlation is not significantly different from zero (e.g.~on the $5\%$ level) and to use $\hat m = \ell^*-1$ as an estimate of $m$. As before, the evaluation of statistically significant deviations from zero must be based on local estimates because potential rate changes can bias the estimates of serial correlations as illustrated in Figure \ref{fig:mestimation} B-D. Panel B shows a simulated spike train according to model (\ref{eq:ss}) with negative first order serial correlations, i.e., $\rho_1<0$, and a rate change point in the middle. The corresponding successive ISIs $\xi_i, \xi_{i+1}$ on which the estimation of $\rho_1$ is based are shown in panel C. The global estimate of $\rho_1$ is not even negative but positive (dashed line in C), whereas the true correlation is indicated by the blue and black lines with negative slope; a phenomenon known as Simpson's paradox. 

We therefore propose to estimate $m$ by splitting up the  process into disjoint sections. In each section, serial correlations up to a maximal  lag are calculated, and systematic deviations from zero are investigated for each lag. These sections should be long enough to provide good estimates of serial correlations, and small enough so that most windows remain unaffected by potential change points. In Figure \ref{fig:mestimation} D, the estimates derived from the local estimators in small windows (black and grey dots) agree well with the true correlation structure (green) of the spike train shown in panel B, whereas the global estimators (red) are highly biased.

\begin{figure}[htbp]
\begin{center}
\includegraphics[width=0.45\textwidth]{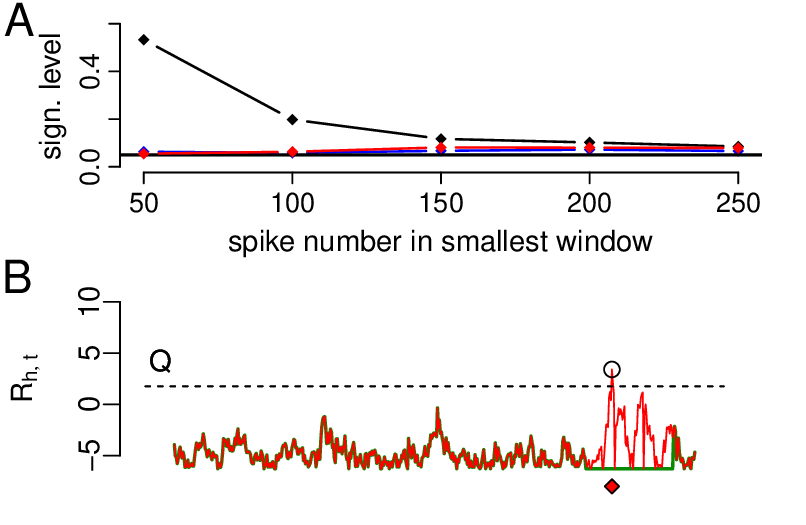}
\caption{Practical applicability to finite data sets and window choice. A.  Empirical significance level of MFT$^{(m)}$ applied to simulations of spike trains from Figure \ref{fig:PSEex}. For positively correlated ISIs (black), the significance level approaches the asymptotic $0.05$ (horizontal line) when increasing the smallest window and thus, the spike number. For negatively correlated ISIs (red, blue), cutting out the $h$-neighborhood of points with undefined $\hat s$, which  occurs particularly in small windows, reduces the number of false positives.  B. lIllustration of test modification. By cutting out the $h$-neighborhood (green curve) around falsely detected change points (red curve) caused in the neighborhood of points with undefined $\hat s$, false positives are reduced for small analysis windows, particularly when $\rho^2/\mu^3$ is close to zero.}
\label{fig:windowsize}
\end{center}
\end{figure}

Finally, we investigate the practical applicability of the proposed procedure to finite windows as it relies on asymptotic thresholds. As mentioned earlier \citep[see also Fig.~$9$ in][]{Messer2014}, simulations suggest that the MFT$^{(0)}$ keeps the asymptotic significance level if the smallest window contains about $100-200$ spikes for spike trains with medium regularity, i.e., if $\sigma^2/\mu^3$ is not too small. If we assume additional covariance structure, we need to consider the term $\rho^2/\mu^3$ instead, which basically determines the asymptotic value of the denominator of $G_{h,t}$. If it takes values close to zero, estimation error may lead to negative estimates of $\rho^2/\mu^3$, in which case  $\hat s$ and $G_{h,t}$ would be not defined. In addition, estimates of $\rho^2/\mu^3$ in the neighborhood may be positive, but extremely small, causing sharp peaks in $G_{h,t}$ and therefore, false positives, particularly when using smaller windows (Fig.~\ref{fig:windowsize} B, red curve with estimated change point). This needs to be taken into account in practice because negative serial correlations may yield very small $\rho^2/\mu^3$. We therefore suggest to slightly modify the MFA by excluding the $h$-neighborhood of points in which the denominator of $G_{h,t}$ is not defined by setting $\hat s:=0$ in this neighborhood (such that $G$ is also set to zero in this case, Fig.~\ref{fig:windowsize} B, green curve). As this has asymptotically no effect, $\hat s$ remains consistent. The empirical significance level of this modified MFT$^{(m)}$ is investigated by application to the three simulated spike trains from Figure \ref{fig:PSEex} by varying the minimal window size. Figure \ref{fig:windowsize} A shows that in these simulations, again about $150-200$ spikes can be required to approximately reach the asymptotic significance level.

\section{Application to spike train recordings}\label{sect:dataana}
We apply the proposed methods, principles and algorithms  to an experimental  data set of spike trains obtained from spontaneous activity recordings of dopaminergic neurons in the substantia nigra and ventral tegmental area of anaesthetized  mice, as described previously \citep{Schiemann2012,Subramaniam2014}. The data set contains $44$ spike trains of length $600$ seconds, with a mean rate of about $4$ spikes per second. The set of analysis windows was therefore chosen as $H=\{50, 75, 100\}$ seconds, yielding an expected number of about $200$ spikes in the smallest window. 

\begin{figure*}[htbp]
\begin{center}
\includegraphics[width=0.21\textwidth]{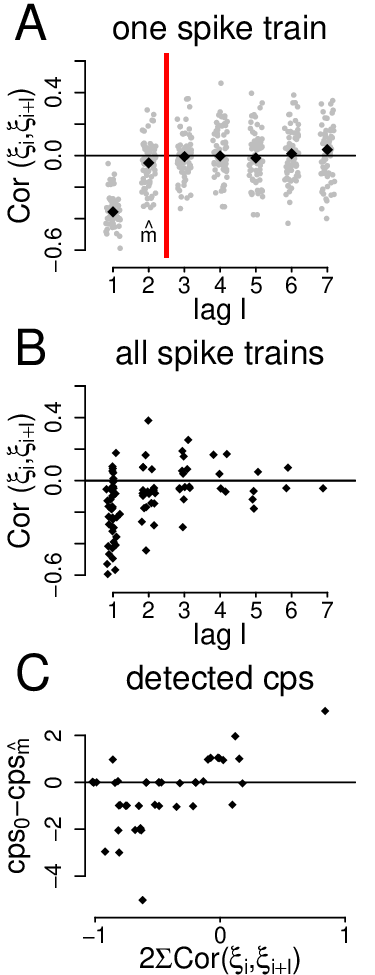}
\includegraphics[width=0.408\textwidth]{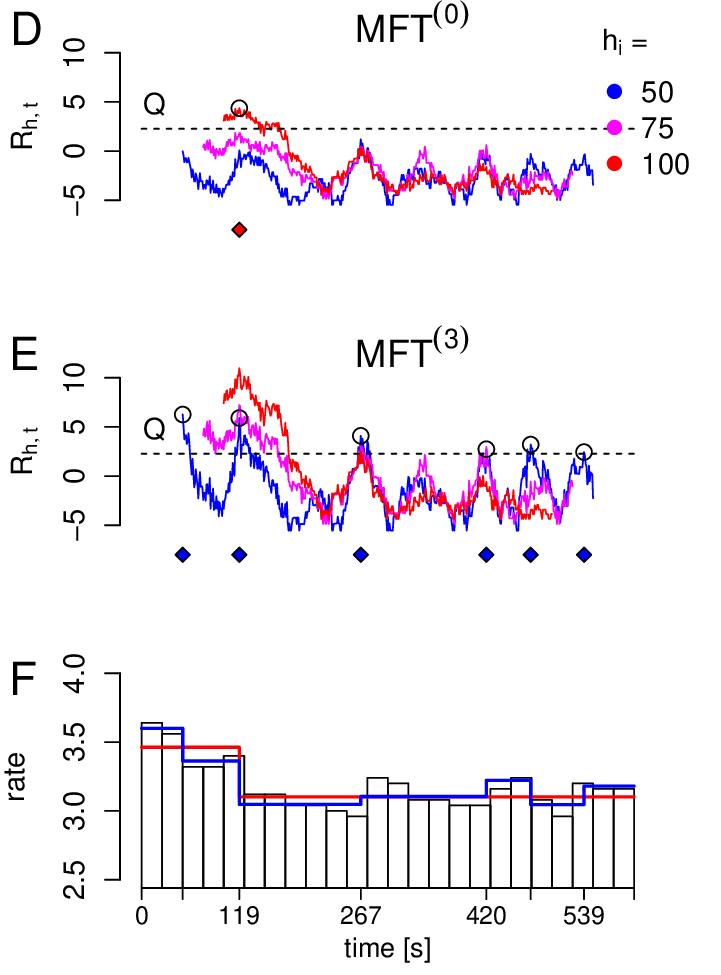}
\caption{Application of the MFA to a data set of spike trains with weak serial correlations. A. Serial correlations estimated in disjoint windows of $50$ life times each (grey), and medians (black). Vertical line indicates cutoff value for $\hat m$ for the respective spike train. B. Serial correlations are short and typically negative. Median serial correlations derived as in A for all spike trains, plotted up to the respective estimate $m_i$ for every spike train $i$. C. Difference between the number of  change points estimated by the  MFA$^{(0)}$ and the  MFA$^{(\hat m)}$, as a function of the contribution of the serial correlations to $\rho^2$ (eq.~(\ref{eq:CovData})). D. Application of the  MFA$^{(0)}$ and E. the  MFA$^{(\hat m)}$ to one spike train with correlation profile similar to A. F. The rate profile of the sample spike train, and the rate profiles estimated by the MFA$^{(0)}$ (red, one estimated change point) and the MFA$^{(\hat m)}$ (blue, six estimated change points).}
\label{fig:dataanalysis}
\end{center}
\end{figure*}

We estimated the maximal lag $\hat m$ for every spike train separately as described in section \ref{sect:application} (Fig.~\ref{fig:mestimation} D). To that end, we used disjoint windows of $50$ ISIs  to estimate serial correlations, and estimated $m+1$ as the first lag for which deviations from zero were not significant on the $5\%$-level using a Wilcoxon test. Figure \ref{fig:dataanalysis} A shows a typical example for one spike train. The serial correlation of lag one shows considerable deviation from zero, the correlation of lag two is  small but still deviating from zero, and all other correlations do not strongly deviate from zero, leading to $\hat m = 2$ for this spike train. The corresponding estimates of serial correlations up to $\hat m_i$ are shown in panel B for all spike trains.  The values of $\hat m$ were $\hat m\le 3$ in about $90$\% of all cases, ranging up to a maximum of $7$, and tended to be negative in the majority of spike trains.

In this more frequent case of negative serial correlations, the MFA$^{(\hat m)}$ typically detected more change points than the MFA$^{(0)}$, leading also to rate profiles matching better with visual inspection (D-F). In order to measure this effect as a function of the degree of serial correlations, we plot the difference between the number of  change points estimated by the MFA$^{(0)}$ and by the MFA$^{(\hat m)}$ in panel C as a function of an estimate of the term
\begin{equation}\label{eq:CovData}
2\sum_{\ell=1}^{m} {\text{Cor}}(\xi_i,\xi_{i+\ell}) =  \frac{ \rho^2-\sigma^2}{\sigma^2},
\end{equation}
which measures the contribution of serial correlations to $\rho^2$. As expected, when this term is negative, the MFA$^{(0)}$ typically estimated much fewer change points, often none at all. In the rare cases where this term is positive, the MFA$^{(0)}$ typically estimated more change points than the  MFA$^{(\hat m)}$.

\section{Discussion}
We have presented a multiple filter test (MFT) that can test the null hypothesis of constant firing rate and estimate change points in the  rate of spike trains especially if these show  dependencies in their ISI structure as is often observed experimentally. Detection of subtle rate changes can be used for extracting meaningful signals from neuronal spike trains and, more generally, it can be an important preprocessing step for statistical analyses that are sensitive to rate changes.  

Our procedure incorporates multiple features that are particularly important for practical application in spike train analysis: (1) un unknown number of rate changes can occur (2) on multiple time scales, (3) other process parameters such as the variance of inter spike intervals can be unknown, and (4) processes can show a high variety of patterns and distributions, including particularly serial dependencies.

The initial version of the MFT for rate change detection introduced in \cite{Messer2014} was developed for renewal processes with a wide range of life time distributions but assumed independence of ISIs, which does often not hold in empirical neuronal spike trains. The MFT uses a filtered derivative process with multiple filters that converges weakly to a parameter free limit process that can be used to obtain the rejection threshold for the test. By specifically estimating serial dependencies in the test statistic, we show that the new MFT can be applied to a variety of empirical firing patterns, including positive and negative serial correlations as well as tonic and bursty firing. Note that  the conditions for the present new MFT include models where the life times are independent or where the life times are dependent but show no serial correlations. In these cases the results of the present MFT would be identical to the results of the original MFT \citep{Messer2014}. This is because zero serial correlation implies that $\rho^2=\sigma^2$ i.e., the terms that are responsible for the difference in the methods are identical.

For practical application, it is necessary to estimate  the denominator of the test statistic, $s$, consistently. We have therefore proposed a consistent local parameter estimator under  $m$-dependence. Although more complex theoretical approaches for consistent estimation are  available for the more general case of ergodicity \citep{Berkes2005,Wu2009,Xiao2012,MuhsalKirch2014}, we focus  on $m$-dependence because it is technically simple and suitable for empirical data analysis with finite spike numbers. Especially under the alternative of rate changes, global estimators of $s$ are affected by rate changes and yield erroneous results. Therefore, our simulations argue strongly for local estimates of $s$ within small windows as these are less affected by potential change points. Even these local estimators require that $m$ is small relative to the window size used for estimation. This implies that even under $m$-dependence the performance can be suboptimal if $m$ is large and change points occur frequently. This is because large $m$ requires large windows with constant firing rate for the estimation of $s$. If change points occur frequently, such windows cannot be found, and consequently, $\hat s$ will be affected by change points within the used estimation windows. Therefore, in practice, only cases with a moderate number of change points and short range dependencies can be considered, i.e., when $m$ is small or serial correlations decay fast with the lag. According to practical examples such as the data set used here \citep[e.g.~Figure \ref{fig:examples} and][]{Ratnam2000,Chacron2001,Nawrot2007,Farkhooi2009}, this is a typical case for empirical spike trains.

One practical limitation of the presented method is its asymptotic nature, which requires a sufficient number of spikes, i.e., about $100-200$ events in the smallest window, which prevents change point detection in shorter time scales. Therefore, it can be considered particularly useful for spontaneous activity, rather than for short trials with many external stimuli or behavioral events. For these cases, different methods such as for example point process adaptive filter methods \citep[e.g.,][]{Eden2004} may be useful. The main problem with using smaller window sizes is that the asymptotic threshold, $Q$, is too low when the smallest window does not contain sufficiently many spikes. One possibility to deal with this issue could be to replace $Q$ by a threshold $Q_b$ derived from a block bootstrap procedure \citep{Singh1981,Goncalves2011,Kreiss2012}, where the block size needs to be chosen such that serial correlations can be treated appropriately. In our simulations of the spike trains used in Figure \ref{fig:windowsize} A, a block bootstrap procedure kept the asymptotic significance level by  increasing the rejection threshold $Q_b$  (data not shown). However, while $Q$ always depends only on  the window set $H$ and the time $T$, under the alternative hypothesis of change points, $Q_b$ largely depends on the properties of the spike train. This can render interpretation difficult in case of change points. Bootstrap can be advantageous when $\rho^2$ is close to zero, for example due to strong negative correlations, such that large amounts of data would need to be excluded from the analysis due to negative estimates of $\hat s$, potentially also including the change points themselves. In such cases, bootstrap procedures can enhance detection probability by avoiding this exclusion. In other cases, detection probability can be reduced, which often makes the use of small windows equally unsatisfactory for bootstrap procedures. In addition, the derivation of $Q_b$ takes considerably longer than the derivation of $Q$. We therefore recommend to use the asymptotic threshold and a minimal spike number of about $100-200$ events in the smallest window, but bootstrap options are also made available in the provided code.

As a second limitation, the present method assumes the rate to be a step function with clear change points. As a consequence, other forms of the rate function, such as ramps or rhythmic behavior, will be described by corresponding step functions.

Our simulations illustrate the necessity of incorporating serial correlation in the MFT. For positive correlations, our new MFT is necessary to reduce the number of false positives, which can be highly enhanced when falsely assuming independence. For the frequent case of negative correlations, these reduce the variability of the spike count and therefore enhance the detection probability of change points, yielding a higher potential of signal extraction from noisy spike trains. Indeed, it has been suggested that sensorial neural systems, such as the electroreceptive organs of weakly electric fish \citep{Chacron2001} and primary somatosensory cortical neurons in rats \citep{Nawrot2007} use this feature to increase their information transfer capacity. In this, our method takes into account a feature of information transfer in point processes with a direct correlate in the actual function of neuronal circuits.

In order to illustrate the  performance of the method, we have applied the new MFA$^{(\hat m)}$ to a data set of empirical spike trains and compared its performance to the classical MFA$^{(0)}$ that falsely assumes independence of ISIs. For all spike trains, serial correlations of small orders were estimated by using small windows to account for potential bias caused by rate changes. In the rare case of positive correlations, the classical MFT$^{(0)}$ that falsely assumes independence detected up to twice as many change points as the new MFT$^{(\hat m)}$. In the more typical case of negative serial correlations, the new MFT$^{(\hat m)}$ detected many more change points than the MFT$^{(0)}$. The new MFT$^{(\hat m)}$ then yielded rate profiles  matching better with visual inspection, indicating a higher detection power of potential neuronal rate signals. Potential applications of our novel algorithm include the extraction of information-rich signals from noisy spike trains, especially when there are no clear behavioral or sensorial triggers, e.g. spontaneous activity recordings. It can also potentially be used as a pre-processing step for other statistical analyses, and for detecting long-term but subtle rate changes, which may reflect transitions of neuromodulatory states \citep{Lee2012}. 

\section*{Acknowledgements}
We would like to thank G\"otz Kersting for helpful comments on weak convergence principles. This work was supported by the German Federal Ministry of Education and Research (BMBF, Funding number: 01ZX1404B) and by the Priority Program 1665 of the German Research Foundation.

\begin{appendix}
	
\section{Proofs}

Here we show consistency of the estimators $\hat s^2$ of $s^2$ in equations (\ref{s1_mdependent}) and (\ref{s2_mdependent}). Recall that these were
\begin{align*}
\text{global estimator:} \quad & 2h\hat\rho^2/\hat\mu^3 \quad & \text{(see Lemma \ref{mdependent1})}\\
\text{local estimator:} \quad & \left(\frac{\hat\rho_{\text{ri}}^2}{\hat\mu_{\text{ri}}^3}+\frac{\hat\rho_{\text{le}}^2}{\hat\mu_{\text{le}}^3}\right)h  \quad & \text{(see Lemma \ref{mdependent2})}
\end{align*}
The used estimators $\hat \rho, \hat \mu, \hat \rho_{\text{le}}, \hat \rho_{\text{ri}}, \hat \mu_{\text{le}}, \hat \mu_{\text{ri}}$ are the empirical means and estimates of $\rho$ given in equations (\ref{eq:rhok}), derived from the whole process in the global estimator and from the local right and left windows at time $t$ in the local estimator.

\begin{lemm}\label{mdependent1}
Let $\{\xi_i\}_{i\ge 1}$ be an $m$-dependent process in $\mathscr P$ and $(\hat s_{nh,nt}^2)_{t}$ the global estimator as in (\ref{s1_mdependent}). Then it holds in $(D[h,T-h],d_{\|\cdot\|})$ almost surely as $n\to\infty$ that 
$$(\hat s_{nh,nt}^2 /n)_{t}\longrightarrow (2h\rho^2/\mu^3)_{t},$$
where $d_{\|\cdot\|}$ denotes the supremum norm. 
\end{lemm}

\noindent\textbf{Proof:} Note that the global estimator $\hat s$ does not depend on $h$ and $t$, i.e., the formulation of $\hat s$ as a process is artificial. We show that $\hat\mu\to\mu$ a.s.~and $\hat\rho_\ell\to\rho_\ell$ a.s.~as $n\to\infty$ for $\ell=0,1,2,\ldots$ where $\rho_0=\sigma^2$. Since $\{\xi_i\}_{i\ge 1}$ is $m$-dependent and square-integrable, the sequence $\{\xi_i\xi_{i+\ell}\}_{i\ge 1}$ is integrable and $(m+\ell)$-dependent, thus ergodic. Then, the ergodic theorem, see e.g., \cite{Klenke2008}, states almost surely as $n\to\infty$
\begin{align}\label{slln1}
\frac{1}{n}\sum_{i=1}^{n} \xi_i \longrightarrow \mathbb{E}[\xi_1]=\mu
\quad\textrm{and}\quad 
\frac{1}{n}\sum_{i=1}^{n} \xi_i\xi_{i+\ell}\longrightarrow \mathbb{E}[\xi_1\xi_{1+\ell}].
\end{align}
Since the life times are a.s.~positive and integrable, it follows $N_{nT}\to\infty$ a.s.~as $n\to\infty$ (cmp.~the proof to Lemma A.1.~in \cite{Messer2014}).  Thus, in $(\ref{slln1})$, the value $n$ can be exchanged with the random number of observations $N_{nT}$ (respectively $N_{nT}-(\ell-1)$).  
Hence, for $n\to\infty$, we find $\hat\mu\to\mu$ a.s.~and $\hat\rho_\ell\to\rho_\ell$ a.s., so that the finite sum $\hat\rho^2\to\rho^2$ a.s. By construction of $\hat s^2$ the statement holds.\hfill$\Box$

\begin{lemm}\label{mdependent2}
Let $\{\xi_i\}_{i\ge 1}$ be an $m$-dependent  process in $\mathscr P$ and for all $T>0$ and $h\in(0,T/2]$ let $((\hat s_{nh,nt})^2)_{t}$ be the local estimator as in (\ref{s2_mdependent}). Then it holds in $(D[h,T-h],d_{\|\cdot\|})$ almost surely as $n\to\infty$ that $((\hat s_{nh,nt})^2 /n)_{t}\to (2h\rho^2/\mu^3)_{t}$. 
\end{lemm}

\noindent\textbf{Proof:} We show the uniform a.s.~convergence of $(\hat{\mu}_{\text{le}})_{t}$ and $(\hat{\mu}_{\text{ri}})_{t}$ to the constant $\mu$ in Lemma \ref{lemm_conv_esti_mu}, and the uniform a.s.~convergence of the summands $(\hat\rho_{\text{le},\ell})_{t}$ and $(\hat\rho_{\text{ri},\ell})_{t}$ of $\hat \rho_{\text{le}}^2$ and $\hat \rho_{\text{ri}}^2$ to the constant $\rho_\ell$ in Lemma \ref{lemm_conv_esti_rho}. This implies the statement, since uniform almost sure convergence interchanges with finite sums in general and with products if the limits are constant.\hfill $\Box$\\

We start with a uniform a.s.~result for the scaled counting process $(N_t)_{t\ge 0}$.
Throughout, we use the following approach: First, we state an almost sure convergence result for the finite dimensional marginals of the processes. This essentially results from the ergodic theorem. Then, by a discretization argument, we show uniform a.s.~convergence.  

\begin{lemm}\label{lemm_conv_nt}
Let  $\{\xi_i\}_{i\ge 1}$ be a process in $\mathscr P$ with $\mathbb{E}[\xi_1]=\mu$.
Then we have in $(D[h,T-h], d_{\|\cdot\|})$
 almost surely as $n\to\infty$ that
\begin{align}
\left(\frac{N_{nt}-N_{n(t-h)}}{nh/\mu}\right)_{t} &\longrightarrow (1)_{t},\label{n000}\\ 
\left(\frac{N_{n(t+h)}-N_{nt}}{nh/\mu}\right)_{t} &\longrightarrow (1)_{t}.\label{conv_nt1}
\end{align} 
\end{lemm}

\noindent\textbf{Proof:}
We show (\ref{conv_nt1}); (\ref{n000}) follows analogously.\\
For  $S_n := \sum_{i=1}^n\xi_i$ for $n\ge 1$, the ergodic theorem implies
 $S_n/ n \to\mu$ a.s. for $n\to\infty$. As we have  $N_t \to \infty$ a.s.~as $t\to\infty$,   $S_{N_t}/N_t\to\mu$ a.s.~as $t\to\infty$. Now, for all $t\ge0$ we find $S_{N_t} \le t \le S_{N_t + 1}$, so that (for all $t$ sufficiently large such that $N_t\ge 1$)
\begin{align*}
\frac{S_{N_t}}{N_t} \le \frac{t}{N_t} \le \frac{S_{N_t+1}}{N_t+1}\frac{N_t+1}{N_t}.
\end{align*}
Since the left hand side and the right hand side tend to $\mu$ almost surely we obtain $N_t / t \to 1/\mu$ a.s.~as $t\to\infty$. For $0\le s<t$, this implies, as $n\to\infty$, almost surely
\begin{align}
\frac{N_{nt} - N_{ns}}{n(t-s)}& = \frac{t}{t-s}\frac{N_{nt}}{nt} - \frac{s}{t-s}\frac{N_{ns}}{ns}\nonumber\\
& \longrightarrow \frac{t}{t-s}\frac{1}{\mu} -\frac{s}{t-s}\frac{1}{\mu} = \frac{1}{\mu}.\label{conv_nt2a}
\end{align}
This implies the convergence of the finite dimensional marginal of (\ref{conv_nt1}).
The uniform convergence follows by a discretization argument analogously to the proof of Lemma A.14 in \cite{Messer2014}.\hfill$\Box$

Next, we show the uniform a.s.~convergence of the estimators  $(\hat\mu_{\text{ri}})_{t}$,  $(\hat\mu_{\text{le}})_{t}$,  $(\hat\sigma_{\text{ri}}^2)_{t}$ and $(\hat\sigma_{\text{le}}^2)_{t}$.

\begin{lemm}\label{lemm_conv_esti_mu}
Let  $\{\xi_i\}_{i\ge 1} \in \mathscr P$ with  $\mu:=\mathbb{E}[\xi_1]$. Then it holds in $(D[h,T-h], d_{\|\cdot\|})$ almost surely as $n\to\infty$ that
\begin{align*}
\left(\hat\mu_{\text{le}}\right)_{t} \longrightarrow (\mu)_{t} \qquad\textrm{and}\qquad (\hat\mu_{\text{ri}})_{t} \longrightarrow (\mu)_{t}.
\end{align*}
\end{lemm}

\noindent\textbf{Proof:}
Again we prove the statement only for the right window.  
We find  $(1/n)\sum_{i=1}^{n}\xi_i\to\mu$ a.s., such that \\
$(1/N_t)\sum_{i=1}^{N_t}\xi_i\to\mu$ a.s. as $n\to\infty$. Then we conclude for all $0< s < t$ (the case $s=0$ being similar) as $n\to\infty$ almost surely
\begin{align}
\frac{1}{N_{nt}-N_{ns}} & \sum_{i=N_{ns}+1}^{N_{nt}} \xi_i \nonumber\\
& = \frac{N_{nt}}{N_{nt}-N_{ns}} \left( \frac{1}{N_{nt}}\sum_{i=1}^{N_{nt}} \xi_i  - \frac{N_{ns}}{N_{nt}}\frac{1}{N_{ns}} \sum_{i=1}^{N_{ns}} \xi_i  \right)\nonumber\\
& \longrightarrow \frac{t}{t-s}\left(\mu - \frac{s}{t} \mu\right) = \mu,\label{conv_hat_mu}
\end{align}
making use of Lemma \ref{lemm_conv_nt}.
Thus, for every fixed $t$ we obtain almost surely as $n\to\infty$
\begin{align}\label{konsistenz_mu_t}
\hat{\mu}_{\text{ri}} =\frac{1}{N_{n(t+h)}-N_{nt}-1}\sum_{i=N_{nt}+2}^{N_{n(t+h)}}\xi_i\longrightarrow \mu.
\end{align}

The a.s.~convergence holds for finitely many $t$ simultaneously.
As above, the uniform convergence follows by a discretization argument analogously to the proof of Lemma A.15 in \cite{Messer2014}.\hfill$\Box$

Now we show the uniform a.s.~convergence of covariance estimators.

\begin{lemm}\label{lemm_conv_esti_rho}
Let  $\{\xi_i\}_{i\ge 1}\in \mathscr P$, and let $\hat\rho_{\text{le},\ell}$ and $\hat\rho_{\text{ri},\ell}$ be the local estimators of $\rho_\ell$ in the left and right window, see (\ref{eq:rhok}),(\ref{s2_mdependent}), for $\ell= 0,1,2\ldots$, where $\rho_0=\sigma^2$. Then in $(D[h,T-h], d_{\|\cdot\|})$ a.s.~as $n\to\infty$ we have
\begin{align*}
\left(\hat\rho_{\text{le},\ell}\right)_{t} \longrightarrow (\rho_\ell)_{t} \qquad\textrm{and}\qquad
(\hat\rho_{\text{ri},\ell})_{t} \longrightarrow (\rho_\ell)_{t}.
\end{align*}
\end{lemm}

\noindent\textbf{Proof:}
Again we conclude $(1/n)\sum_{i=1}^n \xi_i\xi_{i+\ell}\to\mathbb{E}[\xi_1\xi_{1+\ell}]$ a.s. as $n\to\infty$. Using $N_{nT}\to\infty$, we find \\
$(1/N_{nt})\sum_{i=1}^{N_{nt}} \xi_i\xi_{i+\ell}\to\mathbb{E}[\xi_1\xi_{1+\ell}]$ a.s. as $n\to\infty$. With a similar argument as in (\ref{conv_hat_mu}), we find for all $0\le s<t$ almost surely as $n\to\infty$
\begin{align}\label{conv:rho_t}
\frac{1}{N_{nt}-N_{ns}-(\ell+1)}\sum_{i=N_{ns+2}}^{N_{nt}-\ell} \xi_i\xi_{i+\ell}\to\mathbb{E}[\xi_1\xi_{1+\ell}].
\end{align}
Together with the previous Lemma \ref{lemm_conv_esti_mu} this implies the almost sure convergence $\hat\rho_{\text{ri},\ell}\to\mathbb{E}[\xi_1\xi_{1+\ell}]-\mathbb{E}[\xi_1]^2 = \rho_\ell$ for every fixed $t$ and thus for the finite dimensional marginals.

In order to obtain the convergence in $(D[h,T-h],d_{\|\,\cdot\,\|})$, we show a.s.~as $n\to\infty$ that
\begin{align}\label{uniform_xi_i_quad3}
\left(\frac{\mu}{nh}\sum_{i=N_{nt}+2}^{N_{n(t+h)}} \xi_i\xi_{i+\ell}\right)_{t}   \longrightarrow (\mathbb{E}[\xi_1\xi_{1+\ell}])_t.
\end{align}
The convergence of the finite dimensional marginals follows from (\ref{conv:rho_t}) together with Lemma \ref{lemm_conv_nt} and Slutsky's theorem.  
We show the uniform convergence (\ref{uniform_xi_i_quad3}) even for $t\in[0,T-h]$. It suffices to  show almost surely that
\begin{align}\label{uniform_xi_i_quad2}
\lim_{n\to\infty} \sup_{t\in[0,T-h]} \frac{\mu}{nh}\sum_{i=N_{nt}+2}^{N_{n(t+h)}} \xi_i\xi_{i+\ell}   &\le \mathbb{E}[\xi_1\xi_{1+\ell}],\\
\lim_{n\to\infty} \inf_{t\in[0,T-h]} \frac{\mu}{nh}\sum_{i=N_{nt}+2}^{N_{n(t+h)}} \xi_i\xi_{i+\ell} &\ge \mathbb{E}[\xi_1\xi_{1+\ell}].\nonumber
\end{align}
Again, we make use of a discretization argument as in \cite{Messer2014}. We make it explicit here, since the mixing terms $\xi_i\xi_{i+\ell}$ were not explicitly considered in the latter article. 
For an $\varepsilon >0$ with $T/\varepsilon\in\mathbb N$, we decompose the time interval $[0,nT]$ into equidistant sections of length $n\varepsilon$. Using the notation $|\lceil x\rceil|:=\lceil x \rceil +1, x\in \mathbb R,$ we bound
\begin{align*}
& \sup_{t\in[0,T-h]} \frac{\mu}{nh}\sum_{i=N_{nt}+2}^{N_{n(t+h)}} \xi_i\xi_{i+\ell} \\
& \le  \max_{j=0,1,\ldots,T/\varepsilon - |\lceil h/\varepsilon \rceil |} \frac{\mu}{nh}\sum_{i=N_{jn\varepsilon}}^{N_{jn\varepsilon+n|\lceil h/\varepsilon\rceil|\varepsilon }} \xi_i\xi_{i+\ell}\\
&  \le \max_{j=0,1,\ldots,T/\varepsilon - |\lceil h/\varepsilon \rceil |} \frac{\mu}{nh}\sum_{i=N_{jn\varepsilon +nh}}^{N_{jn\varepsilon+n| \lceil h/\varepsilon\rceil | \varepsilon }} \xi_i\xi_{i+\ell}\\
&\qquad + \max_{j=0,1,\ldots,T/\varepsilon - |\lceil h/\varepsilon \rceil |} \frac{\mu}{nh}\sum_{i=N_{jn\varepsilon}}^{N_{jn\varepsilon+nh }} \xi_i\xi_{i+\ell}.
\end{align*}
For any $\delta>0$ we can choose $\varepsilon>0$ so that \\
$ \max_{j=0,\ldots,T/\varepsilon - |\lceil h/\varepsilon \rceil |} (N_{jn\varepsilon + n| \lceil h/\varepsilon \rceil | \varepsilon} - N_{jn\varepsilon+nh})/ (\delta n/\mu) \to 1$ a.s. for $n\to\infty$.
Then, for $n\to\infty$, the first summand in the latter display converges to $(\delta/h)\mathbb{E}[\xi_1\xi_{1+\ell}]$ a.s.~and the second summand to $\mathbb{E}[\xi_1\xi_{1+\ell}]$ a.s., since convergence  (\ref{uniform_xi_i_quad3}) holds for finitely many $t$.
Since $\delta$ can be chosen arbitrarily small, we find the first inequality of (\ref{uniform_xi_i_quad2}). The second one follows analogously.  Thus, the convergence in  $(\ref{uniform_xi_i_quad3})$ follows. We then exchange the normalization according to Lemma \ref{lemm_conv_nt}. Omitting $\ell+1$ summands does not change the limit such that the uniform a.s. convergence of $(\hat\rho_{\text{ri},\ell})_{t}$ is shown.
Analogously, the uniform a.s. convergence of $(\hat\rho_{\text{le},\ell})_{t}$ is shown.\hfill $\Box$

\end{appendix}

\bibliographystyle{apalike}

\bibliography{./paper_dependence_messer}

\end{document}